\documentclass[12pt]{article}
\usepackage{amsmath}
\usepackage{amssymb}
\usepackage{latexsym}

\input{diagrams} 
\diagramstyle[PostScript=dvips,nohug]

\usepackage{array} \newcolumntype{C}{>{$}c<{$}}

\newcommand{\hovedfont}{\normalfont\bfseries}
\usepackage{theorem}
\theorembodyfont{\normalfont\itshape}
\theoremheaderfont{\hovedfont}
   \theoremstyle{plain}

   \theoremstyle{change}
\newtheorem{lemma}{Lemma.}[section]
\newtheorem{prop}[lemma]{Proposition.}
\newtheorem{cor}[lemma]{Corollary.}
   \theorembodyfont{\normalfont}
\newtheorem{eks}[lemma]{Example.}
\newtheorem{BM}[lemma]{Remark.}

\numberwithin{equation}{lemma}

\newtheorem{taller}[lemma]{$\!\!$}
\newenvironment{blanko}[1]%
{\begin{taller}{\hovedfont #1}\normalfont}%
{\end{taller}}
\newenvironment{deff}%
{%
\begin{list}{\em Definition. }%
{\setlength{\labelsep}{0mm}\setlength{\leftmargin}{0mm}%
\setlength{\labelwidth}{0mm}\setlength{\listparindent}{\parindent}%
\setlength{\parsep}{\parskip}\setlength{\partopsep}{0mm}}%
\item%
}%
{%
\end{list}%
}
\newenvironment{bevis}%
{%
\begin{list}{\em Proof. }%
{\setlength{\labelsep}{0mm}\setlength{\leftmargin}{0mm}%
\setlength{\labelwidth}{0mm}\setlength{\listparindent}{\parindent}%
\setlength{\parsep}{\parskip}\setlength{\partopsep}{0mm}}%
\item%
}%
{%
\qed\end{list}%
}
\newenvironment{bevis*}[1]%
{%
\begin{list}{\em #1 }%
{\setlength{\labelsep}{0mm}\setlength{\leftmargin}{0mm}%
\setlength{\labelwidth}{0mm}\setlength{\listparindent}{\parindent}%
\setlength{\parsep}{\parskip}\setlength{\partopsep}{0mm}}%
\item%
}%
{%
\qed\end{list}%
}

\newcommand{\bigfatgreek}[1]{\boldsymbol{#1}} 
\newcommand{\evalclass}{\bigfatgreek\eta} 
\newcommand{\incidenceclass}{\bigfatgreek\iota} 
\newcommand{\tangencyclass}{\bigfatgreek\tau} 
\newcommand{\psiclass}{\bigfatgreek\psi} 
 
\newcommand{\boundary}{\bigfatgreek\beta}
\newcommand{\phiclass}{\bigfatgreek\phi} 

 \renewcommand{\d}{\partial} 
\newcommand{\textprod}[2]{\overset{#2}{\underset{#1}{\textstyle{\prod}}}} 
\newcommand{\XX}{\mathfrak{X}} \renewcommand{\epsilon}{\varepsilon} 
\makeatletter 
\renewcommand{\ldots}{\relax\ifmmode\ldotp\ldotp\ldotp\else$\m@th%
\ldotp\ldotp\ldotp\ $\fi} \makeatother 
\providecommand{\qed}{\hspace*{\fill}\nolinebreak[1]\hspace*{\fill} $\Box$}

\newcommand{\pil} {\rightarrow} \newcommand{\langpil} {\longrightarrow}

\newcommand{\shortsetminus}{\,\raisebox{1pt}{\ensuremath{\mathbb r}\,}} 
\newcommand{\kan}{\mbox{\Large \(\omega\)}} 
\newcommand{\OO}{\ensuremath{{\mathcal O}}} \newcommand{\fy}{\varphi} 
\newcommand{\upperstar}{^{\raisebox{-0.25ex}[0ex][0ex]{\(\ast\)}}} 
\newcommand{\lowerstar}{_{\raisebox{-0.33ex}[-0.5ex][0ex]{\(\ast\)}}} 
\newcommand{\df}{\: {\raisebox{0.255ex}{\normalfont\scriptsize :\!\!}}=} 
\newcommand{\tensor} {\otimes} \newcommand{\strictsubset}{\subsetneq} 
\newcommand{\strictsupset}{\supsetneq} \renewcommand{\dim} 
{\operatorname{dim}}  
 
 \newcommand{\ov}{\overline}

\renewcommand{\P}{\mathbb{P}}  
  
\newcommand{\Q}{\mathbb{Q}}

\newcounter{dummycounter}
\newenvironment{punkt-i}%
{%
   \begin{list}%
   {(\roman{dummycounter})}%
   {\usecounter{dummycounter}%
   \setlength{\itemsep}{0em}\setlength{\parsep}{0em}\setlength{\topsep}{0em}%
   \setlength{\itemindent}{0em}\setlength{\labelwidth}{1.8em}%
   \setlength{\labelsep}{0.6em}\setlength{\leftmargin}{2.4em}}%
}%
{\end{list}}

\begin{document}

\title{
   Recursion for twisted descendants\\
   and characteristic numbers of rational curves
}
\author{
   Joachim Kock\thanks{Supported by the National Science Research Council of 
   Denmark}\\
   {\normalsize Universidade Federal de Pernambuco}\\ 
   {\normalsize Recife, Brazil}%
}
\date{{\normalsize 3rd of February 1999}}%

\maketitle 
\begin{abstract}
   On a space of stable maps, the psi classes are modified by subtracting 
   certain boundary divisors.  The top products of modified psi classes, usual 
   psi classes, and classes pulled back along the evaluation maps are called 
   {\em twisted descendants}; it is shown that in genus 0, they admit a 
   complete recursion and are determined by the Gromov-Witten invariants.  One 
   motivation for this construction is that all characteristic numbers (of 
   rational curves) can be interpreted as twisted descendants; this is 
   explained in the second part, using pointed tangency classes.  As an 
   example, some of Schubert's numbers of twisted cubics are verified.
\end{abstract}

\section*{Introduction}

On the moduli space $\ov M_{g,n}(X,\beta)$ of stable maps, the psi classes are 
the tautological classes $\psiclass_i \df c_1( \sigma_i\upperstar \kan_\pi)$.  
Here $\pi$ denotes the projection from the tautological family and $\sigma_i$ 
is the section corresponding to the $i$-th mark.  The gravitational 
descendants are the top products of the psi classes and the classes pulled 
back along the evaluation maps.  This notion was devised by 
E.~Witten~\cite{Witten} in order to introduce gravity into topological sigma 
models.  Recently there has been a lot of interest in these invariants,
especially stemming from the Virasoro 
conjecture, which states that certain differential operators annihilate the 
generating function of the descendants.  See Getzler~\cite{Getzler:9812} for 
the precise formulation, a survey, and a proof of the conjecture in genus 0.

In enumerative geometry, the descendants have not yet completely found their 
place. Since the psi classes are not invariant under pull-back via forgetful 
maps it is not obvious that they should have enumerative interpretation.
(Note that geometrical conditions are related to the image curve, and thus
should be invariant under oblivion.)

In the present note, with attention restricted to the case $\ov 
M_{0,n}(\P^r,d)$, I propose the notions of {\em modified psi classes}, 
and {\em twisted descendants} as vehicles for bringing the power of descendants
into enumerative geometry. The exposition naturally divides into two parts:

The first part (Sections \ref{prelim}--\ref{recursion}) is on twisted 
descendants: The enumerative defect of a psi class $\psiclass_i$ is identified 
as the sum $\boundary_i$ of all boundary divisors having mark $i$ on a 
contracting twig; let the {\em modified psi class} be $\ov\psiclass_i \df 
\psiclass_i - \boundary_i$.  The {\em twisted descendants} are then the usual 
gravitational descendants twisted by the modified psi classes.  The main 
result (of the first part of the paper) is the recursive formula~\ref{rec-rel} 
for these twisted descendants, sufficient to determine them from the usual 
descendants, and thus from the Gromov-Witten invariants.  The proof follows 
from a couple of results on the boundary (Lemma~\ref{Delta-isom} and 
Lemma~\ref{DeltaDelta}), which I also find interesting in their own right.


The second part (Section~\ref{part2}) is devoted to enumerative geometry.  The 
key point is the use of pointed conditions.  Let $\Phi_i$ be the codimension-2 
class of pointed tangency: for a given hyperplane $H$, it is the locus 
of maps tangent to $H$ exactly at the mark $i$.  It is shown 
(Proposition~\ref{Phi}) that this locus is the zero 
scheme of a regular section of a vector bundle, and its class is
$$
\Phi_i = \evalclass_i( \evalclass_i + \ov\psiclass_i),
$$
where $\evalclass_i \df c_1(\nu_i\upperstar \OO(1))$ denotes the evaluation 
class.  Even though these Phi cycles do not in general intersect properly, it is 
shown (Proposition~\ref{enum}) that their top products have enumerative 
significance: the characteristic numbers (including compound conditions not 
covered by previous techniques) are all sorts of top products of Phi classes 
(with distinct marks) and eta classes --- in other words, they are twisted 
descendants.

As an example I compute some characteristic numbers of twisted cubics.
A few of these numbers date back to Schubert and are verified for the first 
time; others apparently are new.

\bigskip

This work is a by-product of my forth-coming Ph.D.\ thesis; I am indebted to 
my advisor Israel Vainsencher for his support.  Also I would like to thank 
Lars Ernstr\"om for introducing me to gravitational descendants, and for a 
very inspiring stay at Kungliga Tekniska H\"ogskolan in Stockholm. 
Finally, thanks are due to Andreas Gathmann for sending me a copy of his 
thesis; it has been an important source of inspiration.


\section{Preliminaries}\label{prelim}
   
\begin{blanko}{General conventions.} --- Throughout we work over the field of 
complex numbers.  We place ourselves in the space $\ov M_{0,n}(\P^r,d)$ of 
Kontsevich stable maps in genus 0, and we always assume $n>0$.  For 
simplicity, the target space is taken to be $\P^r$; however, all results of 
Sections 1--3 immediately generalise to any homogeneous variety.

If $\mu : C \pil \P^r$ is a stable map, the source curve $C$ is a tree of 
projective lines; we will consequently refer to the irreducible components of 
the source curve as {\em twigs}, reserving the term {\em component} for the 
components of cycles in the moduli space.

Let $\nu_i : \ov M_{0,n}(\P^r,d) \pil \P^r$ denote the evaluation map of the
$i$'th mark, and let $\evalclass_i \df c_1(\nu_i\upperstar \OO(1))$ denote the 
pull-back of the hyperplane class.
\end{blanko}

For ulterior reference, let us collect some standard results on psi classes 
and gravitational descendants.  For generalisations, see Getzler 
\cite{Getzler:9801}, \cite{Getzler:9812}, or Pandharipande~\cite{Pand:9806}.

\begin{blanko}{Tautological families.}
   --- The projection $\pi : \ov M_{0,n+1}(\P^r,d) \langpil \ov 
   M_{0,n}(\P^r,d)$ that forgets the last mark, comes with $n$ canonical 
   sections $\sigma_i$ defined by repeating the $i$'th mark, and stabilising.  
   The image of $\sigma_i$ is the boundary divisor $\Delta_i$ having 
   a twig of degree 0 with just the two marks $i$ and $n+1$.  (Throughout we 
   will abuse of language like this when we mean: the general point of 
   $\Delta_i$ represents a map whose source curve has two twigs, one 
   of which is of degree 0 and carries just the two marks $i$ and $n+1$.)  

   The map $\pi$, together with the new evaluation map $\nu_{n+1} : \ov 
   M_{0,n+1}(\P^r,d) \pil \P^r$, is a tautological family over $\ov 
   M_{0,n}(\P^r,d)$ in the following sense: Given $[\mu]\in\ov 
   M_{0,n}(\P^r,d)$ representing a map $\mu:C \pil \P^r$, the fibre 
   $\pi^{-1}([\mu])$ is a curve canonically isomorphic to $C$, and the 
   sections $\sigma_i$ single out $n$ distinct marks on the fibre, 
   corresponding to the distribution of marks on $C$.  Under this 
   identification, the restriction of $\nu_{n+1}$ to the fibre 
   $\pi^{-1}([\mu])$ is exactly the map $\mu$.  
\end{blanko}

\begin{blanko}{Psi classes.}
   --- Let $\kan_\pi$ be the relative dualising sheaf of $\pi$.  The $i$'th 
   {\em cotangent line} of $\ov M_{0,n}(\P^r,d)$ is the line bundle 
   $\sigma_i\upperstar \kan_\pi$.  The $i$'th {\em psi class} is its first 
   Chern class:
   $$
   \psiclass_i \df c_1(\sigma_i\upperstar \kan_\pi).
   $$
   The psi class $\psiclass_i$ can also be thought of as ``minus'' the normal bundle 
   of $\sigma_i$:
   \begin{equation}\label{normalbundle}
   \psiclass_i = - \sigma_i\upperstar \Delta_i.
   \end{equation}

   The psi class pulls back along the forgetful map like this:
\begin{equation}\label{pull-psi}
   \pi\upperstar \psiclass_i = \psiclass_i - \Delta_i.
\end{equation}

In case $n\geq 3$ there is the following expression for the psi class:
\begin{equation}\label{123}
\psiclass_1 = (1 \mid 2,3)
\end{equation}
where $(1 \mid 2,3)$ denotes the sum of boundary divisors having mark 1 on one 
twig and 2,3 on the other twig.  This formula is a consequence of the 
pull-back formula above, and the obvious fact that the psi classes are trivial 
on $\ov M_{0,3}$.

Finally, a psi class $\psiclass_i$ restricted to a boundary divisor $\Delta$
gives the psi class $\psiclass_i$ of the twig carrying $i$. Precisely,
(see \ref{boundary-setup} for notation):
\begin{equation}\label{restr-psi}
\rho_\Delta\upperstar \psiclass_i = \jmath_\Delta\upperstar \psiclass_i
\end{equation}
where it is understood that the psi class on the right hand side is pulled 
back from $\ov M_A$ (resp.\ $\ov M_B$) if $i \in A$ (resp.\ $i\in B$).
\end{blanko}

\begin{blanko}{Gravitational descendants.}
   --- Top intersections of psi classes and eta classes are the  
{\em gravitational descendants}. The following notation is similar to 
what has become standard:
$$
\langle \ \tau_{u_1}(c_1)\cdots \tau_{u_n}(c_n) \ \rangle_d \df
\int \psiclass_1^{u_1}\evalclass_1^{c_1} \cdots \psiclass_n^{u_n}\evalclass_n^{c_n}
\cap [\ov M_{0,n}(\P^r,d) ].
$$

The pull-back formula~\ref{pull-psi} yields three important equations known as 
the puncture, dilaton and divisor equations.  In each formula, the left hand 
side is on a space with one extra mark.

\noindent Puncture equation:\vspace*{-1ex}
$$
\langle \ \textprod{i=1}{n} \tau_{u_i}(c_i) \cdot \tau_0(0) \ \rangle_d 
=
\sum_{j=1}^n \langle\ \textprod{i=1}{n} \tau_{u_i-\delta_{ij}}(c_i) \ \rangle_d
$$
Dilaton equation:\vspace*{-1ex}
$$
\langle \ \textprod{i=1}{n} \tau_{u_i}(c_i) \cdot \tau_{1}(0) \ \rangle_d
=
(-2+n)\cdot \langle \ \textprod{i=1}{n} \tau_{u_i}(c_i) \ \rangle_d
$$
Divisor equation:
$$
\langle \ \textprod{i=1}{n} \tau_{u_i}(c_i) \cdot \tau_0(1) \ \rangle_d
= 
d \cdot \langle \ \textprod{i=1}{n} \tau_{u_i}(c_i) \ \rangle_d +
\sum_{j=1}^n \langle \ \textprod{i=1}{n} \tau_{u_i-\delta_{ij}}(c_i+\delta_{ij}) \ \rangle_d
$$
In these formulae, $\delta_{ij}$ is the Kronecker delta, and any term with a 
negative exponent on a psi class is defined to be 0.
\end{blanko}

\begin{blanko}{Topological recursion relation.}
   --- Finally, the topological recursion relation applies when there are at 
   least three marks. It is a consequence of the 
   restriction formula~\ref{restr-psi} and the linear equivalence $\psiclass_i 
   = (i \mid j, k )$.  The formula is (assuming $u_1\geq 1$):
\begin{equation}
\langle\ \textprod{i=1}{n} \tau_{u_i}(c_i) \ \rangle_d =
\sum  
\langle\ \tau_{u_1-1}(c_1)\ 
\textprod{\underset{\scriptstyle i\neq 1}{i\in A}}{} \tau_{u_i}(c_i)\ \tau_0(e)
\ \rangle_{d_A} \ \cdot \
\langle\ \tau_0(r\!-\!e) \ 
 \textprod{i\in B}{} \tau_{u_i}(c_i)\ \rangle_{d_B}.
\end{equation}
The sum is over all stable splittings $A \cup B = [n]$,
$d_A + d_B = d$, and over $e=0,\ldots,r$.
\end{blanko}
\begin{blanko}{Reconstruction for descendants.}\label{reconstruction}
   --- Together, these equations are sufficient to reduce any gravitational 
   descendant to a Gromov-Witten invariant: Induction on the number of psi 
   classes: if there three or more marks, use the topological recursion 
   relation to reduce the number of psi classes.  Otherwise, use first the 
   divisor equation (backwards) to introduce more marks.  (Note that this step 
   may introduce rational coefficients.)
\end{blanko}

\section{Boundary results}\label{boundary}

\begin{blanko}{General set-up for boundary divisors.}\label{boundary-setup}
   --- To establish notation, let $\Delta$ be the 
   boundary divisor $(A, d_A \mid B, d_B)$, where $A\cup B = [n]$ and 
   $d_A+d_B=d$.  Each twig corresponds to a moduli space of lower dimension, 
   $\ov M_A \df \ov M_{A \cup x_A}(\P^r,d_A)$ and $\ov M_B \df \ov M_{B \cup 
   x_B}(\P^r,d_B)$, more precisely, $\Delta$ is the image of a birational map 
   $\rho_\Delta : \ov M_A \times_{\P^r} \ov M_B \pil \ov M_{0,n}(\P^r,d)$.  
   The fibred product is over the two evaluations maps $\nu_{x_A} : \ov M_A 
   \pil \P^r$ and $\nu_{x_B} : \ov M_B \pil \P^r$, reflecting the fact that 
   in order to glue, the two maps must agree at the mark.  The fibred product 
   is a subvariety in the cartesian product $\ov M_A \times \ov M_B$; let 
   $\jmath_\Delta$ denote the inclusion.  This set-up and notation is used 
   throughout --- summarised in the following diagram: 
\begin{diagram}[w=11ex,h=5ex,tight]
\ov M_{0,n}(\P^r,d)  &     &     \\
\uTo<{\rho_\Delta \; }     &        &     \\
\ov M_A \times_{\P^r} \ov M_B & \rInto_{\jmath_\Delta}& \ov M_A \times \ov M_B.
\end{diagram}
\end{blanko}


\begin{lemma}\label{Delta-isom}
   --- Let $\Delta = ( A, d_A \mid B, d_B)$ be a boundary divisor such that
   $A$ is non-empty and $d_A < d_B$. Then $\rho_\Delta$ is an isomorphism onto 
   $\Delta$.
\end{lemma}
This lemma completes Lemma 12 of \textsc{FP-notes}~\cite{FP-notes}, where it 
is proved that if both $A$ and $B$ are non-empty then $\rho_\Delta$ is an 
isomorphism onto its image.  Together the two lemmas cover all the cases of 
isomorphism.

Note in particular that whenever a boundary divisor $\Delta$ has a twig of 
degree zero, the map $\rho_\Delta$ is an isomorphism onto $\Delta$.

\begin{bevis}
  --- We need to show that given a map $\mu:C\pil \P^r$ in the divisor 
  $\Delta$, there is only one node (referred to as the official node) that can 
  be the limit of the node of generic elements in $\Delta$; in other words, 
  there is only one way to cut the curve.

   {\sc start:} Due to the lemma of \textsc{FP-notes}, we can suppose 
   $B=\emptyset$.  There is a minimal sub-tree $T_A \subset C$ carrying all 
   the marks.  Let $T_A$ grow:

   {\sc step 1:} Whenever $T_A$ disconnects the rest, then there is only one 
   of the remaining connected component that has total degree $\geq d_B$; let 
   $T_A$ grow, incorporating the others.

   {\sc step 2:} As long as the total degree of $T_A$ is less than $d_A$, we 
   must incorporate the next twig.  If this act disconnects the rest, repeat 
   step 1.

   {\sc stop:} If the total degree of $T_A$ is equal to $d_A$ the growth stops 
   here, and we have found the unique spot to cut the curve.  Indeed, if the 
   next twig $C'$ has positive degree, it is clear that we cannot include it 
   in $T_A$.  If $C'$ has degree zero, then by {\em stability} it must 
   intersect another two connected components, each of which must have 
   positive total degree.  Including $C'$ would disconnect the rest, and 
   consequently force one of the remaining components to be included in $T_A$.  
   Anyway, the degree of $T_A$ would become too high.

   This shows that the map $\rho_\Delta$ is bijective onto its image.  The 
   proof of normality of $\Delta$ is identical to that indicated in {\sc 
   FP-notes} and resorts to the locally rigidified moduli space $\ov 
   M_{0,n}(\P^r,d,\ov t)$: For the divisors lying over $\Delta$ to intersect, 
   there must be a symmetry along some degeneration of the configuration 
   represented by $\Delta$.  The above arguments show such symmetries do not 
   occur, so $\Delta$ is normal.  Since $\rho_\Delta$ is bijective between 
   normal varieties (and we are in characteristic zero) it must be an 
   isomorphism.
\end{bevis}

\begin{lemma}\label{DeltaDelta}
   --- If $\Delta$ is a boundary divisor such that $\rho_\Delta$ is an 
   isomorphism onto its image, then 
   $$
   \rho_\Delta\upperstar \Delta = \jmath_\Delta\upperstar (- \psiclass_{x_A} - \psiclass_{x_B}).
   $$
Here $\psiclass_{x_A}$ (resp.\ $\psiclass_{x_B}$) denotes the pull-back to 
$\ov M_A \times \ov M_B$ from $\ov M_A$ (resp.\ $\ov M_B$) of the psi class of 
the gluing point.
\end{lemma}
This generalises \ref{normalbundle}.  Note in particular that by 
Lemma~\ref{Delta-isom}, the formula holds for all boundary divisors with a 
contracting twig.

For the sake of completeness: If $\Delta$ is not normal, its pull-back along 
$\rho_\Delta$ has one more term which is the push-forth of the singular locus 
of $\Delta$.  I will not prove that formula here, since it is not needed in 
the sequel.
\begin{bevis}
   --- For short, set $\ov M \df \ov M_{0,n}(\P^r,d)$ and $U \df \ov 
   M_{0,n+1}(\P^r,d)$, and let $\pi: U \pil \ov M$ be the forgetful map.  Its 
   restriction $\bar\pi$ to $D_\Delta \df \pi^{-1}(\Delta)$ is clearly a 
   tautological family over $\Delta$; note that $D_\Delta$ is a divisor in $U$ 
   and it has two components $D_A$ and $D_B$, corresponding to whether the new 
   mark is on the $A$-twig or the $B$-twig.  Now the key point is that the 
   family $\bar\pi : D_\Delta \pil \Delta$ has a section $\bar\sigma_x : 
   \Delta \pil D_\Delta$, defined by putting the new mark on the official node 
   and stabilising.  This works because $\Delta$ has only one official node, 
   according to (the proof of) the previous lemma.
   \begin{diagram}[w=6.5ex,h=5.5ex,tight,shortfall=1.3ex]
D_\Delta    & \rInto^{\bar\rho}     & U      \\
\uTo[shortfall=2ex]<{\bar\sigma_x} \dTo>{\bar\pi}     &        & \dTo>{ \ \pi}      \\
\Delta      & \rInto_{\rho}      & \ov M
\end{diagram}
Note that $D_A$ and $D_B$ intersect transversely along $\bar\sigma_x{}\lowerstar 
\Delta$. 

From the cartesian diagram we get
\begin{align}
   \rho\upperstar \rho\lowerstar \Delta & = \rho\upperstar (\pi\lowerstar \bar\rho\lowerstar 
   \sigma_x{}\lowerstar \Delta ) \notag\\
    & = \bar\pi\lowerstar \bar\rho\upperstar (\bar\rho\lowerstar 
                      \sigma_x{}\lowerstar \Delta ) \notag\\
    & = \bar\pi\lowerstar \bar\rho\upperstar (D_A \cap D_B)  \notag\\
    & = \bar\pi_A{}\lowerstar \bar\rho_A\upperstar (D_A \cap D_B) +
    \bar\pi_B{}\lowerstar \bar\rho_B\upperstar (D_A \cap D_B),\label{twoterms}
\end{align}%
where the subscripted maps in the last line are explained by the diagram
\begin{diagram}[w=6ex,h=6ex,tight,tight,shortfall=1.3ex]
      &  & U   &  &  \\
      &\ruInto^{\bar\rho_A}&\uInto>{\bar\rho}&\luInto^{\bar\rho_B}&\\
   D_A& \rInto & D_\Delta& \lInto& D_B\\
      & \rdTo_{\bar\pi_A}&\dTo>{\bar\pi}&\ldTo_{\bar\pi_B}&\\
      && \Delta&&
\end{diagram}
Note that the section $\bar\sigma_x$ factors through $D_A$ giving a section 
$\bar\sigma_{x_A} : \Delta \pil D_A$.  The image of $\bar\sigma_{x_A}$ (which 
we denote by the same symbol) is cut out by the divisor $D_B$ restricted to 
$D_A$.

So, concerning the first term in \ref{twoterms}, we continue
\begin{eqnarray*}
   \bar\pi_A{}\lowerstar ( \bar\rho_A\upperstar (D_A \cap D_B) & = &
   \bar\pi_A{}\lowerstar ( \bar\rho_A\upperstar D_A \cap \bar\sigma_{x_A})\\
    & = & \bar\sigma_{x_A}\upperstar (\bar\rho_A\upperstar D_A )  \\
    & = & \bar\sigma_{x_B}\upperstar (\bar\rho_B\upperstar D_A )  \\
    & = & \bar\sigma_{x_B}\upperstar (\bar\sigma_{x_B}).
\end{eqnarray*}%

By assumption we are identifying $\Delta$ with the fibred product $\ov M_A 
\times_{\P^r}\ov M_B$; under this identification, $\bar\pi_A :D_A \pil 
\Delta$ is the pull-back along $\jmath_\Delta$ of the tautological family $\pi_A : 
U_A \pil \ov M_A$. Similarly for $B$.  The intersection $D_A\cap D_B$ in 
$U$ corresponds to the image of the sections $\sigma_{x_A}$ and $\sigma_{x_B}$, 
along which the pulled-back families $\jmath_\Delta\upperstar U_A$ and 
$\jmath_\Delta\upperstar U_B$ are glued together.
(Since the tautological families do not have the universal property, these 
identifications are only natural up to non-unique isomorphism, but that is 
sufficient for our purposes.)

We find that 
\begin{eqnarray*}
   \bar\sigma_{x_B}\upperstar (\bar\sigma_{x_B}) & = & \jmath_\Delta\upperstar \sigma_{x_B}\upperstar \sigma_{x_B}  \\
    & = & - \jmath_\Delta\upperstar \psiclass_{x_B}.
\end{eqnarray*}%
Same reasoning for the second term in \ref{twoterms} yields the result.
\end{bevis}

\begin{deff}
   --- On a space $\ov M_{0,n}(\P^r,d)$ with $d>0$, let $\boundary_i$ denote 
   the union of all boundary divisors whose general point represents a map 
   contracting the twig carrying the mark $i$. That is, $\boundary_i = 
   \sum_{A\ni i} (A,0 \mid B,d)$.
\end{deff}

\begin{lemma}
   --- Let $\Delta = ( A, 0 \mid B, d)$ be a boundary divisor with a 
   contracting twig.
   
   \noindent
   {\rm (a)} \ If $i \in A$ then
   $$
   \rho_\Delta\upperstar \boundary_i = \jmath_\Delta\upperstar \big(
   \; (i \mid x_A) - \psiclass_{x_A} - \psiclass_{x_B} + \boundary_{x_B}\big).
   $$
   {\rm (b)} \ If $i \in B$ then
   $$
   \rho_\Delta\upperstar \boundary_i = \jmath_\Delta\upperstar \boundary_i.
   $$
\end{lemma}
In (a), the symbol $(i \mid x_A)$ denotes the union of all boundary divisors 
of $\ov M_A$ with mark $i$ on one twig and $x_A$ on the other. The divisor
$\boundary_{x_B}$ lives in $\ov M_B$. In (b), the divisor $\boundary_i$ is 
understood to live in $\ov M_B$, since $i \in B$.
   
\begin{bevis}
   --- Ad.\ (a): Note first of all that $\Delta$ itself is a component in 
   $\boundary_i$.  Therefore the product includes the self-intersection class 
   $\Delta^2 = \jmath\upperstar (-\psiclass_{x_A}-\psiclass_{x_B})$ (by
   \ref{DeltaDelta}).  The 
   other components of $\boundary_i$ are boundary divisors distinct from 
   $\Delta$, and their intersection with $\Delta$ is a reduced effective 
   cycle.  For such a divisor, let $I \subset [n]$ denote the degree 0 part, 
   that is, the part containing $i$.  If neither $I \strictsubset A$ nor $I 
   \strictsupset A$ then clearly the intersection is empty.
   
   Each component of $\boundary_i$ such that $I \strictsubset A$, intersects
   $\Delta$ by breaking the $A$-twig like this:
   
   \vspace*{.5cm}

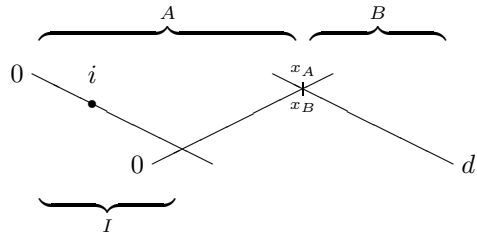
\begin{center}
\setlength{\unitlength}{4mm}
\begin{picture}(16,6)(-16,0)

\put(-1,1.5){\line(-2,1){6}}
\put(-5,4.5){\line(-2,-1){6}}
\put(-9,1.5){\line(-2,1){6}}

\put(-6,3.8){\line(0,1){0.4}}
\put(-6,3.4){\makebox(0,0){{\tiny \(x_B\)}}}
\put(-6,4.6){\makebox(0,0){{\tiny \(x_A\)}}}

\put(-13,3.5){\circle*{.3}}
\put(-13,4.5){\makebox(0,0){{\footnotesize \(i\)}}}

\put(-0.5,1.5){\makebox(0,0){{\footnotesize \(d\)}}}
\put(-11.5,1.5){\makebox(0,0){{\footnotesize \(0\)}}}
\put(-15.5,4.5){\makebox(0,0){{\footnotesize \(0\)}}}

\put(-3.5,6){\makebox(0,0){{\footnotesize \(
   \overset{B}{\overbrace{\phantom{xxxxxxxxx}}}
\)}}} 
\put(-10.5,6){\makebox(0,0){{\footnotesize \(
   \overset{A}{\overbrace{\phantom{xxxxxxxxxxxxxxxxx}}}
\)}}} 
\put(-12.5,0){\makebox(0,0){{\footnotesize \(
   \underset{I}{\underbrace{\phantom{xxxxxxxxx}}}
\)}}} 
\end{picture}  
\end{center}
   such that the mark $i$ is broken apart from the attachment mark $x_A$.  
   There is a 1--1 correspondence between the possible $I \strictsubset A$ and 
   these possible breaks of the $A$-twig.  These codimension-2 boundary cycles 
   correspond again to divisors in $\ov M_A$ such that $i$ is on one twig and 
   $x_A$ is on the other twig.  This accounts for the contribution $(i\mid 
   x_A)$.
   
   Next, the cases $I \strictsupset A$ corresponds to breaking the $B$-twig 
   in this way:
   
   \vspace*{.5cm}

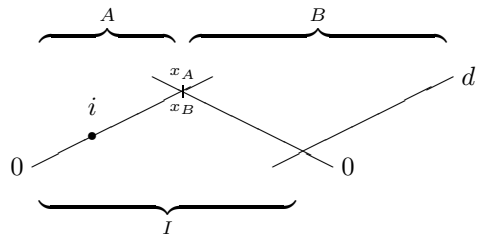
\begin{center}
\setlength{\unitlength}{4mm}
\begin{picture}(16,6)(0,0)

\put(1,1.5){\line(2,1){6}}
\put(5,4.5){\line(2,-1){6}}
\put(9,1.5){\line(2,1){6}}

\put(6,3.8){\line(0,1){0.4}}
\put(6,3.4){\makebox(0,0){{\tiny \(x_B\)}}}
\put(6,4.6){\makebox(0,0){{\tiny \(x_A\)}}}

\put(3,2.5){\circle*{.3}}
\put(3,3.5){\makebox(0,0){{\footnotesize \(i\)}}}

\put(0.5,1.5){\makebox(0,0){{\footnotesize \(0\)}}}
\put(11.5,1.5){\makebox(0,0){{\footnotesize \(0\)}}}
\put(15.5,4.5){\makebox(0,0){{\footnotesize \(d\)}}}

\put(3.5,6){\makebox(0,0){{\footnotesize \(
   \overset{A}{\overbrace{\phantom{xxxxxxxxx}}}
\)}}} 
\put(10.5,6){\makebox(0,0){{\footnotesize \(
   \overset{B}{\overbrace{\phantom{xxxxxxxxxxxxxxxxx}}}
\)}}} 
\put(5.5,0){\makebox(0,0){{\footnotesize \(
   \underset{I}{\underbrace{\phantom{xxxxxxxxxxxxxxxxx}}}
\)}}} 
\end{picture}  
\end{center}
   such that the marks of $I\cap B$ are put on the new middle twig, of degree 0.
   In the space $\ov M_B$ this corresponds to all ways of breaking off a twig
   of degree 0 carrying the mark $x_B$. So these cases give the contribution 
   $\boundary_{x_B}$.
   
   Now for item (b): This time, $\Delta$ has $i$ on the degree $d$ twig, so all
   components in $\boundary_i$ are distinct from $\Delta$. The intersection of
   such a divisor with $\Delta$ is empty unless $I\subset B$ or $I \supset A 
   \cup \{i\}$. 
   
   The cases $I\subset B$ corresponds to breaking the $B$-twig like this
   
   \vspace*{.5cm}

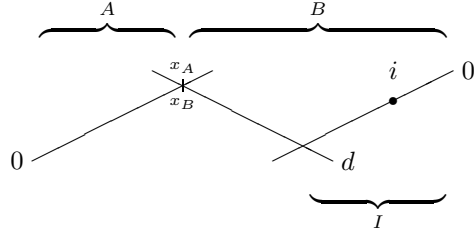
\begin{center}
\setlength{\unitlength}{4mm}
\begin{picture}(16,6)(0,0)

\put(1,1.5){\line(2,1){6}}
\put(5,4.5){\line(2,-1){6}}
\put(9,1.5){\line(2,1){6}}

\put(6,3.8){\line(0,1){0.4}}
\put(6,3.4){\makebox(0,0){{\tiny \(x_B\)}}}
\put(6,4.6){\makebox(0,0){{\tiny \(x_A\)}}}

\put(13,3.5){\circle*{.3}}
\put(13,4.5){\makebox(0,0){{\footnotesize \(i\)}}}

\put(0.5,1.5){\makebox(0,0){{\footnotesize \(0\)}}}
\put(11.5,1.5){\makebox(0,0){{\footnotesize \(d\)}}}
\put(15.5,4.5){\makebox(0,0){{\footnotesize \(0\)}}}

\put(3.5,6){\makebox(0,0){{\footnotesize \(
   \overset{A}{\overbrace{\phantom{xxxxxxxxx}}}
\)}}} 
\put(10.5,6){\makebox(0,0){{\footnotesize \(
   \overset{B}{\overbrace{\phantom{xxxxxxxxxxxxxxxxx}}}
\)}}} 
\put(12.5,0){\makebox(0,0){{\footnotesize \(
   \underset{I}{\underbrace{\phantom{xxxxxxxxx}}}
\)}}} 
\end{picture}  
\end{center}
   --- that is: break off a degree 0 twig containing $i$ but not $x_B$.
   
   The cases $I \supset A 
   \cup \{i\}$ corresponds to breaking off a degree 0 twig containing both $i$ 
   and $x_B$. 
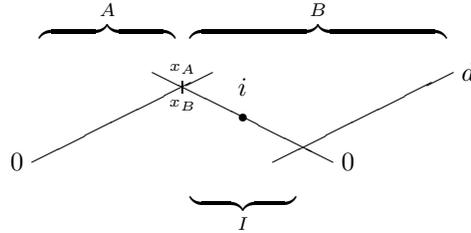
\begin{center}
\setlength{\unitlength}{4mm}
\begin{picture}(16,6)(0,0)

\put(1,1.5){\line(2,1){6}}
\put(5,4.5){\line(2,-1){6}}
\put(9,1.5){\line(2,1){6}}

\put(6,3.8){\line(0,1){0.4}}
\put(6,3.4){\makebox(0,0){{\tiny \(x_B\)}}}
\put(6,4.6){\makebox(0,0){{\tiny \(x_A\)}}}

\put(8,3){\circle*{.3}}
\put(8,4){\makebox(0,0){{\footnotesize \(i\)}}}

\put(0.5,1.5){\makebox(0,0){{\footnotesize \(0\)}}}
\put(11.5,1.5){\makebox(0,0){{\footnotesize \(0\)}}}
\put(15.5,4.5){\makebox(0,0){{\footnotesize \(d\)}}}

\put(3.5,6){\makebox(0,0){{\footnotesize \(
   \overset{A}{\overbrace{\phantom{xxxxxxxxx}}}
\)}}} 
\put(10.5,6){\makebox(0,0){{\footnotesize \(
   \overset{B}{\overbrace{\phantom{xxxxxxxxxxxxxxxxx}}}
\)}}} 
\put(8,0){\makebox(0,0){{\footnotesize \(
   \underset{I}{\underbrace{\phantom{xxxxxxx}}}
\)}}} 
\end{picture}  
\end{center}

   All together we get all the ways of breaking off $i$ with degree 0, that 
   is, the contribution is $\boundary_i$ (from the $B$-space).
\end{bevis}

\section{Modified psi classes}\label{recursion}

\begin{deff}
   --- On a space $\ov M_{0,n}(\P^r,d)$ with $d>0$, the {\em modified psi 
   classes} are by definition
   $$
   \overline \psiclass_i \df \psiclass_i - \boundary_i.
   $$
\end{deff}
The motivation for this definition is the following observation:
\begin{BM}
   --- Let $\pi : \ov M_{0,n+1}(\P^r,d) \pil \ov M_{0,n}(\P^r,d)$ be the 
   universal map that forgets the last mark. Then
   $$
   \pi\upperstar \ov \psiclass_i = \ov \psiclass_i.
   $$
\end{BM}

\begin{BM}\label{mpsi=psi}
   --- On a space with just one mark, we have $\ov \psiclass_1 = \psiclass_1$, 
   so an equivalent definition is $\ov\psiclass_i \df \hat\pi_i\upperstar 
   \psiclass_i$, where $\hat\pi_i$ denotes the map that forgets all marks 
   except $i$.
\end{BM}

\begin{blanko}{Key formula.} \label{restr-mpsi}
   {\em --- Let $\Delta = ( A, 0 \mid B, d)$ be a boundary divisor with a 
   contracting twig. Then
   $$
   \rho_\Delta\upperstar \ov \psiclass_i = 
   \begin{cases}
      \quad \jmath_\Delta\upperstar  \ov \psiclass_{x_B}  \quad& \text{in case }i \in A \\
      \quad\jmath_\Delta\upperstar\ov \psiclass_i \quad &  \text{in case } i \in B.
   \end{cases}
   $$
}\end{blanko}
\begin{bevis}
   --- By definition, $\ov\psiclass_i = \psiclass_i - \boundary_i$.  Suppose 
   $i \in A$. Applying
   the 
   restriction formula~\ref{restr-psi} and the previous lemma
   gives 
   $\sigma_\Delta\upperstar \ov\psiclass_i = \jmath_\Delta\upperstar \big( 
   \psiclass_i  -(i\mid 
   x_A) +\psiclass_{x_A} + \psiclass_{x_B} -\boundary_{x_B} \big)$.  Now note 
   that by stability, the $A$ twig must have yet another mark, say $j$, so 
   by \ref{123} we can write $\psiclass_i=(i \mid j,x_A )$, and   
   $\psiclass_{x_A} = (i,j \mid x_A)$.  It follows that 
   $\psiclass_i + \psiclass_{x_A} = (i\mid x_A)$.  This proves the first case.  
   The second case is immediate from the lemma.
\end{bevis}
The second case shows the geometric nature of the modified classes: If two 
marks come together, a new degree zero twig sprouts out and takes the marks, 
but the corresponding modified psi classes stick to the original twig instead 
of following the marks to the new one, which is illusory from the point of 
view of the image curve.

\subsection{Recursion of twisted descendants}

\begin{deff}
   --- We define a {\em twisted descendant} to be a top product of modified 
   psi classes, usual psi classes, and the eta classes, and agree on the 
   following notation:
   $$
   \langle \ \tau_{u_1}^{m_1}(c_1)\cdots \tau_{u_n}^{m_n}(c_n) \ \rangle_d \df
\int \psiclass_1^{u_1}\ov\psiclass{}_1^{m_1}\evalclass_1^{c_1} \ \cdots \
\psiclass_n^{u_n}\ov\psiclass{}_n^{m_n}\evalclass_n^{c_n}
\ \cap [\ov M_{0,n}(\P^r,d) ].
   $$
\end{deff}
\begin{BM}
   --- Since the modified psi classes are invariant under pull-back along 
   forgetful maps, the puncture, dilaton and divisor equations are also valid
   for twisted descendants: the modified psi classes go through untouched.
   The three equations are not needed in the recursion below.
\end{BM}

\begin{blanko}{Recursion relation.}\label{rec-rel}{\em
   --- Suppose $m_1$, the exponent of the first modified psi class, is positive. 
   Then
$$
\langle \ \textprod{i=1}{n} \tau_{u_i}^{m_i}(c_i) \ \rangle_d =
\langle \ \tau_{u_1+1}^{m_1-1}(c_1) \textprod{i=2}{n}\tau_{u_i}^{m_i}(c_i) \ \rangle_d
- 
\sum \langle \ \textprod{i\in A}{}\tau_{u_i} \cdot \tau_0 \ \rangle
\cdot
\langle \ \tau_{0}^{\mathbf{m}- 1}(\mathbf{c})
\cdot \textprod{i\in B}{}\tau_{u_i}^{m_i}(c_i) \ \rangle_d
$$
The sum is over all partitions $A \cup B = [n]$ with $1\in A$ and $\sharp A 
\geq 2$.   In the last factor, we have 
used the short hand notation $\mathbf{m} \df \sum_{i\in A} m_i$, the sum of 
all the exponents of the modified psi classes corresponding to marks in $A$, 
and similarly we have set $\mathbf{c} \df \sum_{i\in A} c_i$.
}\end{blanko}

The product $\langle \ \textprod{i\in A}{}\tau_{u_i} \cdot \tau_0 \ \rangle$ 
is on the Knudsen-Mumford space $\ov M_{0,A\cup x_A}$ of stable pointed curves 
with labels $A \cup x_A$.  The top products of psi classes on these spaces are 
well known: if $n_A \df \sharp A$ then the factor evaluates to
$$
\frac{(n_A-2)!}{\prod_{i\in A} u_i !}
$$

\begin{bevis}
   --- The left hand side has a factor $\ov \psiclass_1 = \psiclass_1 - 
   \boundary_1$, so we can write 
   it as
   $$
   \psiclass_1\cdot\text{rest } - \; \sum \rho_\Delta\upperstar \text{ rest}
   $$
   where the sum is over all components $\Delta$ of $\boundary_1$.
   Applying the formulae for $\rho_\Delta$, the sum is
   $$
   \sum \int_\Delta \jmath_\Delta\upperstar 
   \bigg( 
   \ov\psiclass_{x_B}^{\mathbf{m}-1} \textprod{i\in B}{} \ov\psiclass{}_i^{m_i} 
   \ \textprod{i=1}{n} \psiclass_i^{u_i}\evalclass_i^{c_i} \big)).   
   $$
   Pulling back via $\jmath_\Delta$ amounts to multiplying with the K\"unneth 
   decomposition of the diagonal (pulled back via the evaluation maps of the 
   gluing marks $x_A$ and $x_B$), and then we can distribute the terms, 
   according to which space they are pulled back from. We get
   $$
\sum
\int_{\ov M_A} \bigg(\textprod{i\in A}{} 
\psiclass_i^{u_i}\evalclass_i^{c_i} \cdot \evalclass_{x_A}^e\bigg) 
\cdot \int_{\ov M_B} \bigg(\ov\psiclass{}_{x_B}^{\mathbf{m}- 1} \evalclass_{x_B}^{r-e}
\cdot \textprod{i \in B}{} \ov\psiclass{}_i^{m_i} \psiclass_i^{u_i} \evalclass_i^{c_i} \bigg)   
   $$
   where the sum is also over $e=0,\ldots,r$.

Finally note that $\ov M_A = \ov M_{0,A\cup x_A}\times \P^r$, so 
the only non-zero contributions come when $e = r - \mathbf{c}$.
\end{bevis}

\begin{cor}\label{cor}
   --- The twisted descendants can be reconstructed from the 
   usual gravitational descendants (and thus from the Gromov-Witten invariants).
\end{cor}
\begin{bevis}
   --- The recursive formula above expresses each twisted descendant as a sum 
   of twisted descendants with either fewer modified psi classes or fewer 
   marks.  When we come down to 1-pointed space, the modified psi classes are 
   equal to the usual psi classes, and then we have a usual gravitational 
   descendant.  (From gravitational descendants to Gromov-Witten invariants
   we have \ref{reconstruction}.)
\end{bevis}

\begin{BM}
   --- Kontsevich and Manin~\cite{KM:9708} also studied certain twisted 
   descendants (called generalised correlators) and obtained a recursion 
   (\cite{KM:9708}, Theorem 1.2) visually simi\-lar to mine.  Instead of 
   twisting with $\ov\psiclass_i$ they twisted with $\phiclass_i\df 
   {\text{st}}\upperstar \psiclass_i$, where $\text{st}$ is the map to 
   Knudsen-Mumford space consisting in forgetting the map and stabilising.  
   Note that $\phiclass_i$ is also a modification of $\psiclass_i$: the 
   difference consists in all boundary divisors such that $i$ is alone on a 
   twig.  
   
   It is interesting to compare the two recursions: The 
   recursion of the present note consists in breaking off twigs containing $i$ 
   and of {\em degree zero}, resulting in a sum over all possible  
   distributions of the remaining marks.  The coefficients are integrals over 
   Knudsen-Mumford spaces.  The recursion of Kontsevich and Manin consists in 
   breaking off twigs containing {\em only} the mark $i$; the resulting sum is 
   over all possible degrees of this twig, and the coefficients are 
   two-pointed descendants.
\end{BM}

\begin{blanko}{Possible generalisations.}\label{speculation}
   --- 
   Let me make an informal comment on possible gene\-ralisations.
   As already noted, all results hitherto remain valid when $\P^r$ is 
   replaced by any homogeneous variety $X$.  Even though the boundary analysis 
   of Section~\ref{boundary} relies on the specific geometry of $\ov 
   M_{0,n}(X,\beta)$ in this case, I believe the results hold also 
   in higher genus and for general smooth projective $X$.  The correction 
   $\boundary_i$ should then be the divisor (in the appropriate component of 
   $\ov M_{g,n}(X,\beta)$) defined as the sum of boundary divisors of type
   \begin{center}\setlength{\unitlength}{5mm}
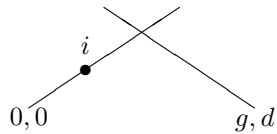
      \begin{picture}(6,3)(0,0)
\put(0,0){\line(3,2){4}}
\put(1.5,1){\circle*{.3}}
   \put(1.5,1.6666){\makebox(0,0){{\footnotesize \(i\)}}}
\put(6,0){\line(-3,2){4}}
   \put(0,-0.3){\makebox(0,0){{\footnotesize \(0,0\)}}}
   \put(6,-0.3){\makebox(0,0){{\footnotesize \(g,d\)}}}
\end{picture}
\end{center}
The correspondingly defined modified psi class should then satisfy the Key 
Formula~\ref{restr-mpsi}, and the Recursion Relation~\ref{rec-rel} would 
follow (with virtual fundamental class instead of the topological one, and 
maybe with signs arising from odd cohomology\ldots).  Note that 
the coefficients $\langle \ \prod \tau_{u_i} \cdot \tau_0 \ 
\rangle$ would still live on the genus 0 space $\ov M_{0,A\cup x_A}$ even in 
the higher genus case.
Such a generalisation would be particularly useful in genus 1 and 2 where 
topological recursion is available (cf.\ Getzler~\cite{Getzler:9801}), in 
which case also the Corollary~\ref{cor} would hold.


\end{blanko}

\section{Tangency conditions and characteristic numbers}
\label{part2}

The motivation for studying twisted descendants comes from enumerative 
geometry:  This final section shows that characteristic 
numbers of rational curves in $\P^r$ are twisted \linebreak descendants.  The 
characteristic numbers are those one can make out of

--- incidence conditions, and

--- conditions of tangency to a hyperplane $H\subset \P^r$ at a 
specified linear subspace of $H$.

\noindent The central idea is the use of pointed conditions.

\bigskip

Previously, working without marks, Pandharipande~\cite{Pand:Q} gave an 
algorithm for computing top intersection products of incidence divisors and 
tangency divisors, allowing for the computation of the simple characteristic 
numbers.  This method comes short in treating compound conditions, since without
marks there is not enough control over the tangencies.
   
The case of $\P^2$ has always been object of particular interest.  Let 
$N_{a,b}$ denote the number of degree $d$ rational plane curves incident to 
$a$ points and tangent to $b$ lines ($a+b=3d-1$).  Pandharipande has shown 
that the corresponding generating function
$$
G(x,y,z)= \sum \sum N_{a,b} y^a x^b \exp(dz)/ a!  b!
$$ 
satisfies the partial differential equation
$$
G_{yz}= -G_x+ G_{xz} - (1/2) G_{zz}^2 + (G_{zz}+ y G_{xz})^2 .
$$
This is the content of the famous unpublished e-mail to Lars 
Ernstr\"om~\cite{Pand:mail}.  The proof is forthcoming in Graber \& 
Pandharipande~\cite{Grab-Pand}.

In another line of development, Ernstr\"om and Kennedy~\cite{EK1}, constructed 
{\em the space of stable lifts} of maps to $\P^2$, and computed characteristic 
numbers of rational plane curves, involving also the codimension-2 condition 
of being tangent to a given line at a specified point.  (It is not clear 
whether the same techniques apply also in higher dimensions.)
      
In a later paper \cite{EK2} they showed how this data is encoded in a ring 
called the {\em contact cohomology ring}, which specialises to the usual 
cohomology ring of $\P^2$ when certain formal parameters are set to zero.
   
They also show how the recursion can be expressed in terms of partial 
differential equations.  It would be very interesting to see if there is a 
relation between my recursion of twisted descendants, and these partial 
differential equations.

\subsection{Pointed tangency condition}

In this section we assume $d\geq 2$ and $r \geq 2$. 
(The propositions~\ref{Phi}, \ref{reduced} and \ref{enum} hold also in the 
case $d=1$, but require slightly different proofs.)

\bigskip

\sloppy

Let $H\subset \P^r$ be a hyperplane, given as the zero scheme of a linear form 
$f\in H^0(\P^r,\OO(1))$.  We consider first the case where there is only one 
mark.  The mark of a map $\mu:C \pil \P^r$ lands in $H$ if and only if 
$\mu\upperstar f$ vanishes at the mark.

\fussy

\begin{deff}
   --- A map $\mu:C\pil\P^r$ with a single mark is {\em tangent} to $H$ at the 
   mark if and only if $\mu\upperstar f$ vanishes at the mark with 
   multiplicity 2, i.e., when also the first derivative vanishes.  Let 
   $\Phi_1(H)$ denote the locus of maps tangent to $H$ at the mark.  (Shortly 
   we will see that the reduced scheme structure is the natural one.)
\end{deff}

In the presence of more marks, we must make the following reservation: in this 
case, the boundary divisor $\boundary_1$ is non-empty, so there exist maps 
whose source curve is reducible with the first mark on a degree $0$ twig, and 
such that the mark (and thus the twig) maps to a point in $H$.  For such a 
map, $\mu\upperstar f$ vanishes identically along the contracting twig and 
thus to any order.  But we do not in general want to consider such 
maps tangent.

Note that these maps have their source collapsed under the forgetful map 
$\hat\pi_1$ that forgets all marks but the first, so a technically convenient way 
to characterise pointed tangency in the case of more marks is this: 
\begin{deff}
   --- A map $\mu : C\pil \P^r$ is tangent to $H$ at the first mark if its 
   image under $\hat\pi_1$ is tangent to $H$ at the mark.  Let $\Phi_1(H)$ be 
   the locus of such maps. Accordingly,
   \begin{equation}
      \Phi_1(H) = \hat\pi_1\upperstar \Phi_1(H),
   \end{equation}
   where the Phi cycle on the right hand side lives on the space with only one 
   mark (compare \ref{mpsi=psi}).
\end{deff}

\begin{prop}\label{Phi}
   --- The locus $\Phi_i(H)$ is the zero scheme of a regular section of a rank 
   two vector bundle $E$, and its class is
   $$
      \Phi_i \df c_2(E) = \evalclass_i (\evalclass_i + \ov \psiclass_i).
   $$
\end{prop}

Since we are going to intersect various Phi cycles, we will need a slightly 
more general set-up.

\begin{deff}
   --- A one-marked map $\mu:C\pil\P^r$ has {\em critical mark} if its 
   differential vanishes at the mark.  The locus of such maps is denoted 
   $\Gamma_1$.  As in the definition of $\Phi_1$, we don't want to consider 
   the maps in $\boundary_1$ as having critical mark, so in a situation with 
   more marks, we define a map to have to have critical first mark if its 
   image under $\hat\pi_1$ has critical mark.  It follows that
   $$
   \Gamma_1 = \hat\pi_1\upperstar \Gamma_1.
   $$
\end{deff}
Note that any map with source of type
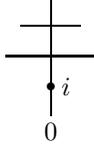
\begin{center}
\setlength{\unitlength}{4mm}
\footnotesize
\begin{picture}(4,5)(0,0)
   \put(2,1){\line(0,1){4}}
   \put(1,4){\line(1,0){2}}
   \put(0.5,3){\line(1,0){3}}
   \put(2,2){\circle*{.3}} \put(2.5,2){\makebox(0,0){\(i\)}}
   \put(2,0.5){\makebox(0,0){\(0\)}}
\end{picture}
\end{center}
has critical $i$'th mark. (Here the 0 indicates the degree of the vertical twig.)

\begin{blanko}{Families of maps.}
   --- Given a flat family of stable maps
\begin{diagram}[w=6ex,h=6ex,tight,shortfall=1.3ex]
\XX      & \rTo^{\mu}      & \P^r      \\
\uTo[shortfall=2ex]<{\sigma_1}\dTo>{\pi}     & \ruTo_{\nu_1}      &     \\
S,    &     & 
\end{diagram}
consider the following conditions:


\begin{quote}\begin{punkt-i}
   \item[{\sf (A)}] $S$ is equidimensional and reduced.
   \item[{\sf (B)}] The locus of maps with reducible source is 
   of pure codimension 1 in $S$ (or it is empty).
   \item[{\sf (C)}] For each mark $i$, the locus $\Gamma_i\subset S$ has 
   codimension at least 2.
\end{punkt-i}\end{quote}
\end{blanko}

Clearly, $\ov M_{0,n}(\P^r,d)$ satisfies all three conditions.

\begin{prop}\label{reduced}
   --- In any family satisfying {\sf (A)}, {\sf (B)} and {\sf (C)}, for 
   general $H$, the scheme $\Phi_1(H)$ is of codimension 2 and satisfies {\sf 
   (A)} and {\sf (B)} (or it may be empty).
\end{prop}

\begin{bevis*}{Proof.}
   --- Let $\hat\pi_1$ denote the map that forgets all marks but the first.
   Since $\Phi_1(H) = \hat\pi_1\upperstar \Phi_1(H)$, and since the properties 
   are preserved under flat pull-back, it is enough to prove the proposition in 
   the case of only one mark.
   
   Consider the pull-back of $f$ along $\mu$, and its first jet relative to 
   $\pi$, (i.e.\ we derive only in vertical direction).  The vanishing of this 
   section
$$
\d_\pi^1\mu\upperstar f \in J_\pi^1\mu\upperstar \OO(1)
$$
defines the locus of maps tangent to $H$.  Pulling back along $\sigma_1$ we 
get only those maps tangent to $H$ at the mark, so set-theoretically 
$\Phi_1(H)$ is the zero locus of the section
$$
\fy \df \sigma_1\upperstar \d_\pi^1(\mu\upperstar f)
$$
of the rank-2 vector bundle
$$
E \df \sigma_1\upperstar J_\pi^1(\mu\upperstar \OO(1)).
$$
(Note that $\pi$ is not a smooth map, so a 
priori the sheaf $J_\pi^1(\mu\upperstar \OO(1))$ needs not be locally free.  
But $\pi$ {\em is} smooth in a neighbourhood of the image of $\sigma_1$; 
therefore $E$ is indeed a vector bundle.)

Off $\Gamma_1$, the vector bundle $E$ is generated by global sections (namely 
the value and the value of the derivative), so by Kleiman's Bertini theorem 
(Remark 5 and 6 in \cite{Kleiman:transversality}), for a general $H$, the 
zero-scheme $Z(\fy)\shortsetminus \Gamma_1$ is reduced of pure codimension 2
(or else it is empty).

Now $\Phi_1(H)$ is contained in the divisor $\evalclass_1(H)$ of maps
whose mark lands in $H$; this divisor moves in a linear system without base 
points, so the intersection $\Gamma_1 \cap \evalclass_1(H)$, and hence
the intersection $\Gamma_1\cap 
Z(\fy)$, is of codimension (at least) 3 in $S$. It follows that the 
whole of $Z(\fy)$ is in fact reduced of pure codimension 2,  so $\Phi_1(H) 
= Z(\fy)$ satisfies {\sf (A)}.


To see that $\Phi_1(H)$ satisfies {\sf (B)}, let $B\subset S$ denote 
the locus of maps with reduced source, and apply the Kleiman-Bertini
to $B\shortsetminus\Gamma_1$. It follows that for general $H$, the intersection 
$\Phi_1(H)\cap B \subset S$ is of codimension 3, away from $\Gamma_1$. But 
anyway $\Phi_1(H)\cap \Gamma_1 \subset S$ is of codimension 3, so the 
dimensionality holds everywhere.
\end{bevis*}
\begin{BM}\label{dim-reduct}
   --- The last dimension reduction argument of the proof above holds for any 
   subvariety $D\subset S$ of codimension at least 1: For general $H$, the 
   intersection $\Phi_1(H) \cap D$ is of codimension at least 3 in $S$.
\end{BM}

\begin{BM}
 --- Note that $\ov M_{0,n}(\P^r,d)$ is Cohen-Macaulay.  Being locally the 
 quotient of a smooth variety by a finite group, its Cohen-Macaulay-ness
 follows from the following algebraic result, valid in characteristic zero: 
 {\em --- If $R$ is a Cohen-Macaulay ring, acted upon by a finite group $G$, 
 then the ring of invariants $R^G$ is also Cohen-Macaulay.} (See Bruns and 
 Herzog~\cite{Bruns-Herzog}, Cor.~6.4.6.)
\end{BM}

\begin{bevis*}{Proof of Proposition~\ref{Phi}.}
   --- Suppose first there is only one mark, and continue the \linebreak 
   notation of the proof above, with $S=\ov M_{0,1}(\P^r,d)$.  Since $Z(\fy)$ 
   is of expected co\-dimension and $\ov M_{0,1}(\P^r,d)$ is Cohen-Macaulay, the 
   section $\fy \in H^0(S,E)$ is regular, and $[Z(\fy)] = c_2(E) \cap [\ov 
   M_{0,1}(\P^r,d)]$.  This top Chern class is easily computed via the 
   fundamental jet bundle exact sequence (pulled back along $\sigma_1$)
   $$
0 \langpil \sigma_1\upperstar \big(\mu\upperstar \OO(1)  \tensor \kan_\pi
\big)
\langpil \sigma_1\upperstar \big(J_\pi^1 \mu\upperstar \OO(1)\big)
\langpil \sigma_1\upperstar \mu\upperstar \OO(1) \langpil 0 ;
$$
the last term is just $\nu_1\upperstar \OO(1)$, with first Chern 
class $\evalclass_1$, and the first term is also a line bundle, with first 
Chern class $\evalclass_1 + \psiclass_1$. Since $\psiclass_1 = \ov\psiclass_1$ 
on a space with just one mark, this proves the proposition in case $n=1$.

In the case of more marks, take the vector bundle to be the pull-back of $E$
along the forgetful map $\hat\pi_1$. Since the classes $\Phi_1$, 
$\ov\psiclass_1$ and $\evalclass_1$ are preserved under this pull-back, the 
result follows.
\end{bevis*}
\newcommand{\fd}{={\raisebox{0.255ex}{\normalfont\scriptsize\!\!:}}\: }
\begin{BM}
   --- Invoking the Cohen-Macaulay-ness of $\ov M_{0,n}(\P^r,d)$ is not 
   strictly necessary to obtain the class of $\Phi_1$; in fact the formula 
   $\Phi_1 = \evalclass_1(\evalclass_1+\ov\psiclass_1)$ holds in any family 
   satisfying {\sf (A)}, {\sf (B)} and {\sf (C)}.  To see this, recall (from 
   Fulton~\cite{Fulton}, 14.1) that there is a \linebreak localised top Chern 
   class ${\bf Z}(\fy) \in A\lowerstar (Z(\fy))$ whose image in $A\lowerstar 
   (S)$ is $c_2(E) \cap [S]$; when $Z(\fy)$ has expected codimension, with 
   cycle $[Z(\fy)] \fd \sum m_i [D_i]$, then ${\bf Z}(\fy) = \sum e_i [D_i]$, 
   and $1\leq e_i \leq m_i$ (\cite{Fulton}, Example 14.4.4).  In our case, 
   $Z(\fy)$ is of correct dimension and is furthermore {\em reduced}, hence 
   $[Z(\fy)] = {\bf Z}(\fy) = c_2(E) \cap [S]$.
\end{BM}

\subsection{Characteristic numbers}

We will show that the top intersections of Phi classes (with distinct marks) 
together with eta classes have enumerative significance. The key result is
that for general choices of hyperplanes $H$, the corresponding top 
scheme theoretic intersection is reduced and (at most) zero dimensional. 

However, the transversality argument has an interesting extra twist due to the 
\linebreak phenomenon of {\em illusory excess}: For any $3\leq t \leq r$, the 
intersection cycle $\Phi_1\cdots \Phi_t$ (which ought to be of codimension 
$2t$) has excess along the locus of maps with source
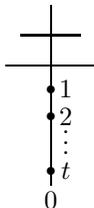
\begin{center}
\setlength{\unitlength}{4mm}
\footnotesize
\begin{picture}(4,7)(0,-2)
   \put(2,-1){\line(0,1){6}}
   \put(1,4){\line(1,0){2}}
   \put(0.5,3){\line(1,0){3}}
   \put(2,2.2){\circle*{.3}} \put(2.5,2.2){\makebox(0,0){\(1\)}}
   \put(2,1.3){\circle*{.3}} \put(2.5,1.3){\makebox(0,0){\(2\)}}
   \put(2.5,0.7){\makebox(0,0){\(\vdots\)}}
   \put(2,-0.5){\circle*{.3}} \put(2.5,-0.5){\makebox(0,0){\(t\)}}
   \put(2,-1.5){\makebox(0,0){\(0\)}}
\end{picture}
\end{center}
where the vertical twig has degree 0 and maps to the intersection of the $t$ 
hyperplanes.  This locus has codimension only $t+2$.  However this is illusory 
excess in the sense that the excessive freedom is just the freedom of marks 
moving on a contracting twig.  The corresponding image curves move much 
less.  Indeed, we can forget $t-1$ marks without destabilising, so in the 
space with no marks, the bad locus gets codimension $t+1$, while the image of 
the cycle $\Phi_1\cdots \Phi_t$ has codimension $t$.  If the total 
intersection is of top degree, a dimension reduction argument then shows that 
the bad curves do not occur at all!  The proof of the Lemma below shows 
how to cope with the phenomenon.

\begin{prop}{Enumerative significance.}\label{enum}
   --- Let $Z\subset \ov M_{0,n}(\P^r,d)$ denote the scheme theoretic 
   intersection of Phi loci and eta divisors, such that for each mark there is 
   at most one Phi, and such that the total codimension of the cycles 
   equals $\dim \ov M_{0,n}(\P^r,d)$.  Then for general hyperplanes $H$, the 
   intersection $Z$ is reduced of dimension 0 (or it may be empty).
   
   Furthermore, the corresponding finitely many maps are immersions with 
   irreducible source (so in particular they have no automorphisms), and they 
   are only simply tangent to the given hyperplanes.
\end{prop}

The proof will follow easily from the following result:
\begin{lemma}
   --- In a $2t$-parameter family satisfying {\sf (A)}, {\sf (B)} and {\sf 
   (C)}, for general hyperplanes $H_i$, the scheme-theoretic intersection 
   $\Phi_1(H_1) \cap \cdots \cap \Phi_t(H_t)$ is reduced of dimension 0 (or it 
   may be empty).
\end{lemma}
\begin{bevis*}{Proof of the Lemma.}
   --- Induction on $t$.  If $t=1$ then the result follows from 
   Proposition~\ref{reduced}.  Suppose the lemma holds for families of 
   dimension $2t-2$ and consider a family $S$ of dimension $2t$.
   For short, set $\Phi_i = \Phi_i(H_i)$.
   
   Again by the Proposition~\ref{reduced}, the variety $\Phi_t$, with the 
   induced family, has dimension $2t-2$ and satisfies {\sf (A)} and {\sf (B)}.  
   However, it will not in general satisfy {\sf (C)}, because the locus $P$
	of maps with source
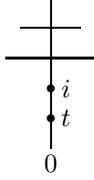
\begin{center}
\setlength{\unitlength}{4mm}
\footnotesize
\begin{picture}(4,6)(0,-1)
   \put(2,0){\line(0,1){5}}
   \put(1,4){\line(1,0){2}}
   \put(0.5,3){\line(1,0){3}}
   \put(2,2){\circle*{.3}} \put(2.5,2){\makebox(0,0){\(i\)}}
   \put(2,1){\circle*{.3}} \put(2.5,1){\makebox(0,0){\(t\)}}
   \put(2,-0.5){\makebox(0,0){\(0\)}}
\end{picture}
\end{center}
	and with the vertical contracting twig mapping to $H_t$, has codimension 3 
	in $S$.  These maps are tangent to $H_t$ and have critical $i$'th mark.

   But we can circumvent this obstacle by forgetting the $t$'th mark: Since
   $P$  has at least four special points on the contracting 
   twig, its image under $\pi_t$ retains the codimension 3 with respect to 
   the whole space, and therefore it has codimension 2 with respect to 
   $\pi_t(\Phi_t)$.  
   \begin{diagram}[w=6.5ex,h=5ex,tight,shortfall=1.3ex]
P     & \overset{\text{\tiny codim 1}}{\subset}    & \Phi_t& \overset{\text{\tiny codim 2}}{\subset} & S    \\
\dTo     &        & \dTo      \\
\pi_t(P)    & \overset{\text{\tiny codim 2}}{\subset}    & \pi_t(\Phi_t)
\end{diagram}
   So $\pi_t(\Phi_t)$ satisfies {\sf (C)}, and since $\Phi_t \pil 
   \pi_t(\Phi_t)$ is generically one-to-one, the image also satisfies 
   {\sf (A)} and {\sf (B)}.
   
   By induction, for general choices of hyperplanes $H_1,\ldots,H_{t-1}$, 
   the intersection \linebreak $\Phi_1\cap \cdots\cap \Phi_{t-1}$ in the space 
   $\pi_t(\Phi_t)$ is reduced of dimension zero, and by the dimension 
   \linebreak reduction argument~\ref{dim-reduct} all these finitely many points 
   correspond to simple tangencies.  Now in an open subset around these 
   points, the map $\Phi_t \pil \pi_t(\Phi_t)$ is immersive and one-to-one, so 
   we conclude that also the intersection upstairs $\Phi_1\cap \cdots \cap 
   \Phi_{t-1} \cap \Phi_t$ in the space $S$ is reduced.
\end{bevis*}
\begin{bevis*}{Proof of Proposition~\ref{enum}.}
   --- The space $\ov M_{0,n}(\P^r,d)$ satisfies the assumptions 
   {\sf (A)}, {\sf (B)} and {\sf (C)}.  The intersection of 
   all the eta divisors will also satisfy the three conditions, since the 
   divisors move in base point free linear systems.  We have arrived at a 
   $2t$-parameter family with $t$ Phi conditions, so the lemma applies.
   
   As to the further properties, they follow from the dimension reduction 
   argument: In each step, the loci of maps with reducible source, critical
   marks, or non-simple tangency (flexes or bitangents) have codimension
   at least 1, so when we come down to dimension zero they no longer occur.
\end{bevis*}

Since the intersection is reduced of dimension zero, the number of points in 
it is equal to the top product of the corresponding classes.  The enumerative 
meaning is this.  Each mark $i$ corresponds to a condition; the factor $\Phi_i 
\evalclass_i^{c_i}$ corresponds to the condition of being tangent to a given 
hyperplane $H_i$ at a specified linear subspace of codimension $c_i$ in $H_i$.  
This understood, the Proposition implies that
\begin{cor}
   --- The top products of Phi classes (with distinct marks) and eta classes
   count what they are expected to.
\end{cor}
By Proposition~\ref{Phi}, these numbers are top products of modified 
psi classes and eta classes, so they are twisted descendants, and can be 
computed by the recursion described in \ref{cor}.  I have implemented the 
algorithm in Maple --- the code is available upon request.

\begin{BM}(Relation with the unmarked conditions.)
   --- Let $\incidenceclass$ be the incidence divisor (of maps incident to a 
   given codimension-2 plane), and let $\tangencyclass$ 
   be the tangency divisor (of maps tangent to a given hyperplane), as defined by 
   Pandharipande~\cite{Pand:Q}.
   
   Let $\pi_i : \ov M_{0,n+1}(\P^r,d) \pil \ov M_{0,n}(\P^r,d)$ be the map that 
   forgets mark $i$. Then
   \begin{eqnarray*}
      \pi_i{}\lowerstar \Phi_i & = & \tangencyclass  \\
      \pi_i{}\lowerstar (\evalclass_i^2) & = & \incidenceclass.
   \end{eqnarray*}%
Therefore, for $a+b= rd+r+d-3$ we have
$$
\int \evalclass_1^2 \cdots \evalclass_a^2 \cdot \Phi_{a+1}\cdots \Phi_{a+b} \cap 
[\ov M_{0,a+b}(\P^r,d)]
=
\int \incidenceclass^a  \tangencyclass^b \cap 
[\ov M_{0,0}(\P^r,d)],
$$
by the projection formula.
\end{BM}

\begin{eks}(Simple characteristic numbers of rational plane curves.)  
   --- For $a+b= 3d-1$, let $N_d(a,b)$ denote the number of rational plane 
   curves passing through $a$ points and tangent to $b$ lines.  By the 
   corollary, it is the number
$$
\int \evalclass_1^2 \cdots \evalclass_a^2 \cdot \Phi_{a+1}\cdots \Phi_{a+b} \cap 
[\ov M_{0,a+b}(\P^2,d)].
$$
\end{eks}

\begin{eks}
   --- In $\P^2$, the condition of being tangent to a given line at a 
   specified point is simply $\Phi_i \evalclass_i$.  As in Ernstr\"om and 
   Kennedy~\cite{EK1}, let $N_d(a,b,c)$ denote the number of rational curves 
   passing through $a$ points, tangent to $b$ lines, and tangent to $c$ lines 
   at specified points.
   
   By the corollary, the characteristic number $N_d(a,b,c)$ is the top 
   intersection
   $$
   \int_{\ov M_{0,a+b+c}(\P^2,d)} \evalclass_1^2\cdots \evalclass_a^2
   \cdot \Phi_{a+1}\cdots\Phi_{a+b} \cdot \Phi_{a+b+1}\evalclass_{a+b+1} \cdot
   \Phi_{a+b+c}\evalclass_{a+b+c}
   $$ 
\end{eks}

Computing these top products by recursion of twisted descendants 
reproduces the numbers of Ernstr\"om and 
Kennedy~\cite{EK1}.

\subsection{Twisted cubics}

Just to close the exposition with some {\em numbers}, I present three tables 
of characteristic numbers of twisted cubics in $\P^3$, most of which, to the 
best of my knowledge, could not have been computed with previous techniques.  
The discussion above immediately yields the computability of the 
characteristic numbers $N_3(a,b,c,d,e)$ one can make out of the following five 
conditions: ($a$) incident to a line; ($b$) tangent to a plane; ($c$) passing 
through a point; ($d$) tangent to a plane at a specified line; and ($e$) 
tangent to a plane at a specified point.  There are 249 such numbers; 7 of 
which are Gromov-Witten invariants; 42 others have $d=e=0$.  Most of these 49 
numbers were computed by Schubert~\cite{Schubert}, and they were verified to 
modern standards of rigour by Kleiman, Str\o mme, and Xamb\'o~\cite{KSX}.

Although Schubert considers also the two compound conditions ($d>0$ or 
$e>0$), he only computes 6 out of the remaining 200 numbers, namely those with 
$c=4$ or $c=5$, and he attributes them to Sturm:
   \begin{eqnarray*}
      N_3(0,0,5,1,0) & = & 2  \\
      N_3(1,0,4,0,1) & = & 3  \\
      N_3(0,1,4,0,1) & = & 6  \\
      N_3(2,0,4,1,0) & = & 17  \\
      N_3(1,1,4,1,0) & = & 34  \\
      N_3(0,2,4,1,0) & = & 68,
   \end{eqnarray*}%
The techniques of the present note confirm these numbers;  for example, the number 
$N_3(1,1,4,1,0)$ of twisted cubics passing 4 points, incident to 1 line, 
tangent to 1 plane, and tangent to 1 plane at a given line is
   $$
   \evalclass_1^3\evalclass_2^3\evalclass_3^3\evalclass_4^3\cdot
   \evalclass_5^2 \cdot \evalclass_6(\evalclass_6+\ov\psiclass_6)
   \cdot \evalclass_7^2(\evalclass_7+\ov\psiclass_7) = 34.
   $$
   
The number $N_3(1,0,4,0,1) = 3$ was previously confirmed by 
Gathmann~\cite{Gathmann} using Gromov-Witten invariants of projective space 
blown up in a point.  (His method yields, for any $r\geq 2$, $d\geq 2$, the 
number of rational curves passing one specified point with tangent direction 
contained in a specified linear space, and incident to further linear spaces 
to get correct codimension.)
   
\begin{eks}
    {\em --- Conditions w.r.t.\ planes and points:} the numbers $N_3(0,b,c,0,e)$
    of twisted cubics
   
   --- tangent to $e$ planes at specified points,
   
   --- passing through $c$ points, and
   
   --- tangent to further $b=12-3e-2c$ planes (at any point):
   
\vspace*{1.5ex}

\begin{center}
\begin{tabular}{C|C|C|C|C|C|} 
                   & e=0    & e=1   & e=2   & e=3   & e=4   \\ 
\hline         c=0 & 56960  & 5552  & 816   &  64   & 4  \\ 
\cline{1-6}    c=1 & 25344  & 3024  & 240   &  12   &\multicolumn{1}{c}{} \\ 
\cline{1-5}    c=2 & 11968  &  888  &  40   &\multicolumn{2}{c}{} \\ 
\cline{1-4}    c=3 &  3376  &  136  &  2    & \multicolumn{2}{c}{} \\ 
\cline{1-4}    c=4 &   480  &    6  &\multicolumn{3}{c}{}\\ 
\cline{1-3}    c=5 &    20  & \multicolumn{4}{c}{}\\ 
\cline{1-2}    c=6 &     1  & \multicolumn{4}{C}{b=12-3e-2c}\\ 
\cline{1-2}
\end{tabular}
\end{center}

\vspace*{1.5ex}

\noindent The numbers of the first column are well-known through 
\cite{Schubert} and \cite{KSX}, and could also be computed using the algorithm 
of \cite{Pand:Q}.
\end{eks}

\vspace*{0.5ex}

\begin{eks}
       {\em --- Conditions w.r.t.\ planes and lines:} the numbers $N_3(a,b,0,d,0)$
    of twisted cubics

   --- tangent to $d$ planes at specified lines,
   
   --- tangent to further $b$ planes (at any point), and
   
   --- incident to $a= 12 - 2d- b$ lines:
 
\vspace*{1.5ex}

\begin{center}
\begin{tabular}{C|C|C|C|C|C|C|C|} 
                  &  d=0   &  d=1   &  d=2   &  d=3   &  d=4   &  d=5   &  d=6   \\ 
\hline       b=0  &  80160 &  38568 &  17472 &  7234  &  2630  &  805   &  217   \\ 
\cline{1-8}  b=1  & 134400 &  60208 &  24720 &  8972  &  2780  &  730   &\multicolumn{1}{c}{} \\ 
\cline{1-7}  b=2  & 209760 &  85216 &  30864 &  9616  &  2552  &  580   &\multicolumn{1}{c}{} \\ 
\cline{1-7}  b=3  & 297280 & 107024 &  33472 &  8952  &  2064  & \multicolumn{2}{c}{} \\ 
\cline{1-6}  b=4  & 375296 & 117312 &  31584 &  7344  &  1504  &\multicolumn{2}{c}{}\\ 
\cline{1-6}  b=5  & 415360 & 112128 &  26272 &  5424  &\multicolumn{3}{c}{}\\ 
\cline{1-5}  b=6  & 401920 &  94528 &  19680 &  3712  &\multicolumn{3}{c}{}\\ 
\cline{1-5}  b=7  & 343360 &  71776 &  13664 &\multicolumn{4}{c}{} \\ 
\cline{1-4}  b=8  & 264320 &  50464 &   9152 &\multicolumn{4}{c}{}\\ 
\cline{1-4}  b=9  & 188256 &  34032 &\multicolumn{5}{C}{a= 12 - 2d- b}\\ 
\cline{1-3}  b=10 & 128160 &  22688 &\multicolumn{5}{c}{}\\ 
\cline{1-3}  b=11 &  85440 &\multicolumn{6}{C}{}\\ 
\cline{1-2}  b=12 &  56960 &\multicolumn{6}{c}{}\\ 
\cline{1-2}
\end{tabular}
\end{center}

\vspace*{1.5ex}

\noindent Again, the first column is well-known through \cite{Schubert}, 
\cite{KSX} and \cite{Pand:Q}.
\end{eks}


\begin{eks}
       {\em --- Combinations of the three tangency conditions:} the numbers 
       \\ $N_3(0,b,0,d,e)$
    of twisted cubics
   
   --- tangent to $e$ planes at specified points,
   
   --- tangent to $d$ planes at specified lines, and
   
   --- tangent to further $b=12-3e-2d$ planes (at any point):
   
\vspace*{1.5ex}

\begin{center}
\begin{tabular}{C|C|C|C|C|C|} 
                     & e=0    & e=1    & e=2    & e=3    & e=4    \\ 
\hline         d=0   & 56960  &  5552  & 816    &  64    &  4     \\ 
\cline{1-6}    d=1   & 22688  &  2320  & 240    &  16    &\multicolumn{1}{c}{}\\ 
\cline{1-5}    d=2   & 9152   &  904   &  72    &\multicolumn{2}{c}{}\\ 
\cline{1-4}    d=3   & 3712   &  328   &  20    &\multicolumn{2}{c}{}\\ 
\cline{1-4}    d=4   & 1504   &  110   &\multicolumn{3}{c}{}\\ 
\cline{1-3}    d=5   &  580   & \multicolumn{4}{c}{}\\ 
\cline{1-2}    d=6   &  217   & \multicolumn{4}{C}{b=12-3e-2d}\\ 
\cline{1-2}
\end{tabular}
\end{center}

\vspace*{1.5ex}

\noindent (The first number in this table is well-known\ldots)
\end{eks}


\medskip

\noindent{\sc Joachim Kock\\
Departamento de Matem\'atica\\
Universidade Federal de Pernambuco\\
Cidade Universit\'aria, 50670-901\\
Recife --- PE --- Brasil}\\[2mm]
e-mail: jojo@dmat.ufpe.br

\end{document}